\documentclass{article}
\usepackage[english,russian]{babel}

\usepackage{amssymb}
\usepackage{amsmath}

\usepackage[dvips]{graphicx}

\begin{document}


\begin{center}
 \textbf{Investigation of the solutions of the Cauchy problem and boundary-value problems for the ordinary differential equations with continuously changing order of the derivative}
\end{center}

\begin{center}
 \textbf{N.A. Aliyev${}^{1}$, R.G. Ahmadov${}^{1}$}

\noindent ${}^{1}$Baku State University, Baku, Azerbaijan

\noindent e-mail: ahmadov\_ramiz@hotmail.com
\end{center}

\noindent 

\noindent Abstract. As is known, the problems for the differential equations with continuously changing order of the derivatives are not considered completely. In this paper we consider the initial and boundary value problems for this type of linear ordinary differential equations with constant coefficients and obtain analytic representations of solutions of these problems.

It should be noted that this area is one of the less studied fields of modern mathematics and there are not effective methods for the study of problems for such differential equations, just as we study the problem for partial differential equations, with both with additive and a multiplicative derivatives.

The method used here is based on an invariant for the above mentioned derivative, i.e. on the functions that do not change for the derivative of any real positive order.

\noindent \textbf{Keywords.} Cauchy problem, boundary value problem, ordinary differential equations, constant coefficients, changing order of the derivative, initial conditions, manifolds.

\noindent \textbf{AMS Subject Classification:   }

\textbf{ }

\section{Introduction}

In the additive analysis both in discrete [1], [2] and continuous cases [3], [4], various problems are well studied [5], [6]. Not so extensive, but still investigated as discrete [7], [8] and continuous problems for the differential equations with multiplicative derivatives [9], [10]. For the partial differential equations, the number of initial conditions under consideration coincides with the order of the derivative of the equation with respect to time [11]. Concerning the boundary conditions for partial differential equations, the number of non-local boundary conditions (if the boundary of the domain is divided into two parts) also coincides with the order of the derivative with respect to the space variables [12].

If we consider the ordinary differential equations of fractional order, then the number of conditions not related to the order of the derivative in this equation, but to the step of change of the order of the derivative [13].

We here consider as the Cauchy problem and boundary value problems for the linear ordinary differential equations with constant coefficients, where the order of the derivative is changed continuously.  In solution of this problem we use mainly invariant function technique. 

Euler was the first to construct an invariant with respect to usual additive derivative function $e^{x} $. Then he constructed the function $e^{\lambda x} $by the help of which ordinary linear homogeneous differential equations with constant coefficients he corresponded some algebraic equation (characteristic equation), i.e. he indeed algebrized the differential equation [14].

Then, with the help of Mittag-Leffler function the invariant function for fractional derivative was constructed [15], [16].

In this paper, using the Volterra function [17] the invariant is constructed for the derivative, the order of which is changing continuously.

\section{Equation and its general solution. }

\noindent 

\noindent Consider the equation 
\begin{equation} \label{GrindEQ__1_} 
\int _{0}^{\sqrt{2} }D^{\alpha } y(x)d\alpha =0, x>a>0, 
\end{equation} 
where $\sqrt{2} $ is an order of equation,  the order of which changes continuously from zero to $\sqrt{2} $.

\noindent Following the Volterra function [17], consider the function that by $\lambda =1$ is an invariant for any real order 

\noindent 
\begin{equation} \label{GrindEQ__2_} 
h_{+0} (x,\lambda )=\int _{-1}^{\infty }\frac{x^{\upsilon } }{\upsilon !}  \lambda ^{\upsilon } d\upsilon , 
\end{equation} 
where $\upsilon !$ by $\upsilon >-1$ is defined by the  Euler's gamma function $\Gamma $[15] and is assumed that  $\upsilon !=\infty $ if 
\begin{equation} \label{GrindEQ__3_} 
\upsilon \le -1. 
\end{equation} 

Note that the function \eqref{GrindEQ__2_} by $\lambda =1$ is a Volterra function if the integral takes the beginning not from -1 and from zero. Now it is easy to see that for any${}_{0<\alpha \in R}$${}_{,}$
\begin{equation} \label{GrindEQ__4_} 
D^{\alpha } h_{+0} (x,\lambda )=\lambda ^{\alpha } h_{+0} (x,\lambda ).                                                  
\end{equation} 
Thus the solution of the equation \eqref{GrindEQ__1_} we seek as 
\begin{equation} \label{GrindEQ__5_} 
y(x)=h_{+0} (x,\lambda )\, . 
\end{equation} 
Substituting \eqref{GrindEQ__5_} into \eqref{GrindEQ__1_} we get
\[\int _{0}^{\sqrt{2} }\lambda ^{\alpha } h_{+0} (x,\lambda ) d\alpha =0,\] 
or
\begin{equation} \label{GrindEQ__6_} 
\int _{0}^{\sqrt{2} }\lambda ^{\alpha }  d\alpha =0.                                                               
\end{equation} 

Note that the equation \eqref{GrindEQ__6_} is a generalization of the polynomial when the order of the variable changes continuously. 

It easy to see that calculation of the integral for the characteristic equation \eqref{GrindEQ__6_} leads us to 
\begin{equation} \label{GrindEQ__7_} 
\frac{\lambda ^{\sqrt{2} } -1}{\ln \lambda } =0, 
\end{equation} 
i.e.
\begin{equation} \label{GrindEQ__8_} 
\lambda _{k} =e^{\sqrt{2} k\pi i} ,      k\in Z,\; \; k\ne 0.                                                                    
\end{equation} 

Thus considering \eqref{GrindEQ__8_} from \eqref{GrindEQ__5_} for the general solution of the equation \eqref{GrindEQ__1_} we get
\begin{equation} \label{GrindEQ__9_} 
y(x)=\sum _{\mathop{k\in Z}\limits_{k\ne 0} }c_{k} h_{+0} (x,\lambda _{k} ) , 
\end{equation} 
where $c_{k} $are arbitrary constants and 
\begin{equation} \label{GrindEQ__10_} 
h_{+0} (x,\lambda _{k} )=\int _{-1}^{\infty }\frac{x^{\upsilon } }{\upsilon !}  \lambda _{k} ^{\begin{array}{l} {} \\ {\upsilon } \end{array}} d\upsilon \, \, ,   k\in Z,\; \; k\ne 0,                          
\end{equation} 
are partial solutions of the equation  \eqref{GrindEQ__1_}.

\noindent 

\section{Cauchy problem}

Now let us consider Cauchy problem for the equation \eqref{GrindEQ__1_} with the following initial conditions  
\begin{equation} \label{GrindEQ__11_} 
D^{\alpha } y(x)\left|_{x=a} \right. =\varphi (\alpha ),  \alpha \in {\rm [0,}\sqrt{{\rm 2}} {\rm ).} 
\end{equation} 

To define the constants $c_{k} $ included to the general solution of the equation \eqref{GrindEQ__9_} we substitute its general solution into \eqref{GrindEQ__11_}. Then we obtain 
\begin{equation} \label{GrindEQ__12_} 
\sum _{\mathop{k\in Z}\limits_{k\ne 0} }c_{k} e^{^{\sqrt{2} k\pi \alpha i} } h_{+0} (a,\lambda _{k} ) =\varphi (\alpha ),   \alpha \in [0,\sqrt{2} ].              
\end{equation} 
The functions $e^{\sqrt{2} k\pi \alpha i} $  form orthogonal system by $\alpha \in [0,\sqrt{2} ]$ i.e. 
\begin{equation} \label{GrindEQ__13_} 
\begin{array}{l} {\int _{0}^{\sqrt{2} }e^{\sqrt{2} k\pi \alpha i} e^{\sqrt{2} n\pi \alpha i}  d\alpha =\int _{0}^{\sqrt{2} }e^{\sqrt{2} (k+n)\pi \alpha i}  d\alpha =\frac{e^{\sqrt{2} (k+n)\pi \alpha i} }{\sqrt{2} (k+n)\pi i} \left|\begin{array}{l} {\sqrt{2} } \\ {0} \end{array}\right. =} \\ {\frac{e^{2(k+n)\pi i} -1}{\sqrt{2} (k+n)\pi i} =\left\{\begin{array}{l} {0,\; \; k+n\ne 0,} \\ {\sqrt{2} ,\; k+n=0.} \end{array}\right. } \end{array} 
\end{equation} 
Multiplying \eqref{GrindEQ__12_} by ${}_{e^{-\sqrt{2} n\pi \alpha i} }$ and integrating from zero ${}_{\sqrt{2} }$ we find 
\[\sum _{\mathop{k\in Z}\limits_{k\ne 0} }c_{k} h_{+0} (a,\lambda _{k} ) \int _{0}^{\sqrt{2} } e^{\sqrt{2} k\pi \alpha i} e^{-\sqrt{2} n\pi \alpha i} d\alpha =\int _{0}^{\sqrt{2} }\varphi (\alpha )e^{-\sqrt{2} n\pi \alpha i}  d\alpha ,\] 
or
\[\sum _{\mathop{k\in Z}\limits_{\begin{array}{l} {k\ne 0} \\ {k\ne n} \end{array}} }c_{k} h_{+0} (a,\lambda _{k} ) \frac{e^{2(k-n)\pi i} -1}{\sqrt{2} (k-n)\pi i} +c_{n} h_{+0} (a,\lambda _{n} )\sqrt{2} =\int _{0}^{\sqrt{2} }\varphi (\alpha )e^{-\sqrt{2} n\pi \alpha i}  d\alpha .\] 
Taking into account \eqref{GrindEQ__13_} one may get
\begin{equation} \label{GrindEQ__14_} 
c_{n} =\frac{1}{\sqrt{2} h_{+0} (a,\lambda _{n} )} \int _{0}^{\sqrt{2} }\varphi (\alpha )e^{-\sqrt{2} n\pi \alpha i}  d\alpha ,   n\in Z,   n\ne 0.           
\end{equation} 
Substituting \eqref{GrindEQ__14_} into \eqref{GrindEQ__9_} we define the solution of the Cauchy problem \eqref{GrindEQ__1_}, \eqref{GrindEQ__11_} in the form
\begin{equation} \label{GrindEQ__15_} 
y(x)=\sum _{\mathop{k\in Z}\limits_{k\ne 0} }\frac{h_{+0} (x,\lambda _{k} )}{\sqrt{2} h_{+0} (a,\lambda _{k} )}  \int _{0}^{\sqrt{2} }\varphi (\alpha )e^{-\sqrt{2} k\pi \alpha i}  d\alpha .                       
\end{equation} 
Thus it was proved 

\noindent \textbf{Theorem.} Let $\varphi (\alpha )$ be continuous function. Then the problem   \eqref{GrindEQ__1_}, \eqref{GrindEQ__11_} has a solution that may be presented in the form \eqref{GrindEQ__15_}.

\section{Boundary problem}

\noindent 

Here we consider the equation \eqref{GrindEQ__1_} on ${}_{(a,b)}$ with boundary condition
\begin{equation} \label{GrindEQ__16_} 
a_{0} D^{\alpha } y(x)\left|_{x=a} \right. +b_{0} D^{\alpha } y(x)\left|_{x=b} \right. =\varphi (\alpha ), \alpha \in [0,\sqrt{2} ]{\rm \; ,\; \; \; \; \; \; } 
\end{equation} 
where $\sqrt{2} $ is an order of the differential equation \eqref{GrindEQ__1_}, $a_{0} ,b_{0} $ are constants, $\varphi (\alpha )$is given real valued smooth function. 

In this case the arbitrary constants $c_{k} $ included into the general solution of the equation \eqref{GrindEQ__9_} are determined by substituting \eqref{GrindEQ__9_} into the boundary condition \eqref{GrindEQ__16_}. Then we have
\[a_{0} \sum _{\mathop{k\in Z}\limits_{k\ne 0} }c_{k}  D^{\alpha } y_{k} (x)\left|_{x=a} \right. +b_{0} \sum _{\mathop{k\in Z}\limits_{k\ne 0} }c_{k}  D^{\alpha } y_{k} (x)\left|_{x=b} \right. =\varphi (\alpha ),\] 
or
\begin{equation} \label{GrindEQ__17_} 
\sum _{\mathop{k\in Z}\limits_{k\ne 0} }c_{k}  \lambda _{k}^{\alpha } [a_{0} h_{+0} (a,\lambda _{k} )+b_{0} h_{+0} (b,\lambda _{k} )]=\varphi (\alpha .)                               
\end{equation} 
Following \eqref{GrindEQ__8_} let us consider the functions 
\begin{equation} \label{GrindEQ__18_} 
\lambda ^{\alpha } _{k} =e^{\sqrt{2} k\pi \alpha i} ,      k\in Z,\; \; k\ne 0,    \alpha \in [0,\sqrt{2} ){\rm \; } 
\end{equation} 
It is easy to see that by $k\ne m$
\[(\lambda _{k}^{\alpha } ,\lambda _{m}^{\alpha } )=\int _{0}^{\sqrt{2} }e^{\sqrt{2} k\pi \alpha i} e^{-\sqrt{2} n\pi \alpha i}  d\alpha =\int _{0}^{\sqrt{2} }e^{\sqrt{2} (k-m)\pi \alpha i}  d\alpha =\frac{e^{^{\sqrt{2} (k-m)\pi \alpha i} } }{\sqrt{2} (k-m)\pi i} \left|_{\alpha =0}^{\sqrt{2} } =\right. \] 
\[=\frac{e^{^{2(k-m)\pi \alpha i} } -1}{\sqrt{2} (k-m)\pi i} =0,\] 
and
\[(\lambda _{m}^{\alpha } ,\lambda _{m}^{\alpha } )=\sqrt{2} ,\] 
i.e. the functions \eqref{GrindEQ__18_} are orthogonal.

\noindent Then going back to \eqref{GrindEQ__17_} we get
\[\sum _{\mathop{k\in Z}\limits_{k\ne 0} }c_{k}  [a_{0} h_{+0} (a,\lambda _{k} )+b_{0} h_{+0} (b,\lambda _{k} )]\int _{0}^{\sqrt{2} }\lambda _{k}^{\alpha } \lambda _{-m}^{\alpha } d\alpha =\int _{0}^{\sqrt{2} }\varphi (\alpha )\lambda _{-m}^{\alpha } d\alpha ,  \] 
or
\begin{equation} \label{GrindEQ__19_} 
c_{m} \sqrt{2} [a_{0} h_{+0} (a,\lambda _{m} )+b_{0} h_{+0} (b,\lambda _{m} )]=\int _{0}^{\sqrt{2} }\varphi (\alpha )\lambda _{-m}^{\alpha } d\alpha \, \, .  
\end{equation} 
If
\begin{equation} \label{GrindEQ__20_} 
a_{0} h_{+0} (a,\lambda _{m} )+b_{0} h_{+0} (b,\lambda _{m} )\ne 0,m\in Z,\; \; m\ne 0,                           
\end{equation} 
Then from \eqref{GrindEQ__19_} we obtain
\begin{equation} \label{GrindEQ__21_} 
c_{m} =\frac{\int _{0}^{\sqrt{2} }\varphi (\alpha )\lambda _{-m}^{\alpha } d\alpha , }{\sqrt{2} [a_{0} h_{+0} (a,\lambda _{m} )+b_{0} h_{+0} (b,\lambda _{m} )]} ,m\in Z,\; \; m\ne 0.                      
\end{equation} 
Putting computed in \eqref{GrindEQ__21_} constants ${c_{m} }$ into \eqref{GrindEQ__9_} we obtain the solution of the problem   \eqref{GrindEQ__1_},\eqref{GrindEQ__16_} in the form
\[y(x)=\sum _{\mathop{k\in Z}\limits_{k\ne 0} }\frac{\int _{0}^{\sqrt{2} }\varphi (\alpha )\lambda _{-k}^{\alpha } d\alpha , }{\sqrt{2} [a_{0} h_{+0} (a,\lambda _{k} )+b_{0} h_{+0} (b,\lambda _{k} )]} h_{+0} (x,\lambda _{k} ) ,\quad x\in [a,b].\quad (22)\] 
Finally for the considered boundary problem \eqref{GrindEQ__1_}, \eqref{GrindEQ__16_} the following theorem is proved. 

\noindent \textbf{Theorem.} If $0<\alpha <b,\quad a_{0} ,\; b_{0} $are given numbers, $\varphi (\alpha )$is continuous real function and the relation function \eqref{GrindEQ__20_} is valid, then solution of the boundary problem \eqref{GrindEQ__1_},\eqref{GrindEQ__16_} exists and may be presented in the form of \eqref{GrindEQ__22_}. \textbf{}

\noindent \textbf{Note. }The case when the coefficients $a_{0} ,\; b_{0} $of the boundary condition \eqref{GrindEQ__16_} are the functions of the variable $\alpha $ is an open problem.

\end{document}